\documentclass[11pt]{article}
\usepackage{proof}
\usepackage{bussproofs}
\usepackage{amssymb}
\usepackage{amsmath}
\usepackage{booktabs}
\usepackage{enumerate} 
\usepackage{xcolor}
\usepackage{setspace}
\usepackage{pstricks}
\usepackage{url}
\usepackage{varwidth}
\usepackage[left=3cm,top=3cm,bottom=3cm,right=3cm]{geometry}
\usepackage{amsthm}
\usepackage{amssymb}
\theoremstyle{definition}

\title{Hierarchical Multiverse of Sets}
\author{Ahmet \c{C}evik}
\date{}
\begin{document}
\maketitle

\begin{abstract}
In this paper, I develop a novel version of the multiverse theory of sets called {\em hierarchical pluralism} by introducing the notion of `degrees of intentionality' of theories. The presented view is articulated for the purpose of reconciling epistemological realism and the multiverse theory of sets so as to preserve a considerable amount of epistemic objectivity when working with the multiverse theory. I give some arguments in favour of a hierarchical picture of the multiverse in which theories or models are thought to be ordered with respect to their plausibility, as a manifestation of endorsing the idea that some set theories are more plausible than others. The proposed multiverse account settles the {\em pluralist's dilemma}, the dichotomy that there is a trade-off between the richness of mathematical ontology and the objectivity of mathematical truth. The view also extends and serves as an alternative position to Balaguer's {\em intention-based Platonism} from which he claims that a certain version of mathematical pluralism follows.
\end{abstract}


\vspace{0.5cm}

I will begin with a short review of set-theoretical pluralism along with some contemporary accounts. I will then argue why having a hierarchical picture of theories is a plausible view and that it is possible to have a reasonable realist based pluralism, that realism in truth value and pluralism in set theory can be fairly consistent with each other.\footnote{The term pluralism in this paper will generally refer to set-theoretic multiversism.} After summarising how the pluralist philosophy of mathematics emerged from geometry, I discuss how pluralism manifested itself in set theory as a result of the existence of independent statements. I will move on to the multiverse theory of sets and the programmes dedicated for its study. I will then elaborate more on hierarchical pluralism to reconcile the multiverse theory of sets and epistemological realism, i.e., {\em realism in truth value} in Shapiro's terms. For this I offer an order theoretical solution to the problem what I call the {\em pluralist's dilemma}, the dichotomy that there is a trade-off between the richness of mathematical ontology and objectivity of mathematical truth. An attempt to reconcile realism and pluralism was also pursued by Mark Balaguer in his \cite{BalaguerPluralism}, where he claims that his {\em intention-based Platonism} entails a version of mathematical relativism. The work presented here can therefore be considered as an alternative view to Balaguer's position in relating realism to pluralism. However, I will also discuss what hierarchical pluralism provides and how it differs from Balaguer's account.

\section{Plurality of Models and the Multiverse Theory}

Mathematical pluralism, historically speaking, began with questions about the conception of `space' in different contexts. The discovery of non-Euclidean geometries in the early 19th century undermined the belief that the Euclidean perception of geometry was the unique model associated with the notion of `space'. The key notion behind {\em pluralism} is the change of `context'. Our first experience of mathematical pluralism perhaps stems from the interpretation of Euclid's fifth postulate of geometry in different contexts of `space'. 
The early 19th century was a major turn in geometry for understanding the mathematical conception of space. It was shown independently by Janos Bolyai, Nikolai Lobachevsky, and Carl Friedrich Gauss that self-consistent non-Euclidean spaces could be defined in which the parallel postulate does not hold. 
Even so that the parallel postulate holds in the Euclidean geometry, it does not hold in non-Euclidean geometries. The parallel postulate is true in some models, yet false in other models. In other words, the parallel postulate is {\em independent} from the other four axioms of Euclidean geometry. The existence of non-Euclidean geometries supports the {\em pluralist} conception of mathematics (or geometry to be more specific), whereas the {\em monist} conception---the view that only a single absolute model should be adopted---argues for that the truth or falsity of mathematical propositions are justified within the frame of one absolute model.

Naturally the first question one might wonder about the pluralist conception of geometry is then which geometry is actually true. Henri Poincar\'{e} \cite{Poincare1902} says that `one geometry cannot be more true than another; it can only be more convenient'.\footnote{Poincar\'{e}, p. 50, 1902.} Knowing that there are many different geometries, which geometry should be adopted? Should we choose, in the first place, an absolute and intuitively plausible model of geometry? We usually tend to pick our axioms based on how useful they are to the mathematical discipline. This approach finds its roots in naturalism.\footnote{See Maddy \cite{Maddy1997} (1997) for a detailed account on mathematical naturalism. For Maddy's multiverse concerns, see her \cite{Maddy2017}, also see Ternullo \cite{Ternullo2019} (2019).} Naturalism suggests that the universe we pick should be determined pragmatically and be based on the fruitfulness of the universe.\footnote{There is a strong connection between pluralism and the question whether or not we should add more axioms in mathematics. A theory, i.e., a set of sentences, may have more than one models. The pluralist is not obliged to choose a particular model. But do all models need to have equal `status'? From the viewpoint of some known studies and what I will present here, the answer is no.} These are expected to rely on what naturalists call {\em extrinsic justifications}. If we were to rely on {\em intrinsic justifications}, i.e., properties that are intuitively plausible, in order to decide which model to choose, then we would be more inclined to endorse a monist view, thinking that our justifications would have to be based on a standard meaningful structure which complies with our mathematical intuition. There can certainly be two intrinsically justified yet contradictory statements, Axiom of Choice (AC) and the Axiom of Determinacy (AD) for instance. However, the view that will be presented here is based on the idea that their `degrees of intentionality' are comparable to each other, that our comparisons about two competing theories can be in favour of one of them. The reason for this is that every intrinsically justifiable statement is supported along with extrinsic features, and we can be sure which theory is preferrable over the other in case we cannot determine this solely by intrinsic qualities of the theory. For example, even if we assume AC and AD have the same intrinsically appealing qualities, since we have more experience with ZFC than ZF+AD in practice, we may say that, in this state of affairs, ZFC is more preferrable than ZF+AD.

\subsection{Set-theoretic arena}

The pluralist idea that there may exist different interpretations of a mathematical concept is not restricted to just geometry, as we can also generalise the pluralist view to sets, hence have a multiverse theory of sets. Although pluralist interpretations of set theory go as far back as to Mirimanoff \cite{Mirimanoff1917} and his studies on non-well-founded sets, results appeared in the later decades have given rise to pluralism in a broader sense in the set-theoretic arena. The ground breaking achievement was the independence of the Continuum Hypothesis (CH) from the axioms of ZFC set theory, which was proved by Kurt G\"{o}del \cite{GodelCH} and Paul Cohen \cite{Cohen}. There are many different forms of stating CH. But one common way of putting it is as follows. 
\vspace{0.5cm}

\noindent {\em Continuum Hypothesis:} Any subset of real numbers is either countable or has the same cardinality as that of the set of real numbers. 
\vspace{0.5cm}


The standard interpretation of the set concept---under the standard interpretation of the set theoretic language---is generally thought to be based on the iterative conception of sets, where every set is ranked in the so-called {\em cumulative hierarchy} with respect to its order type. The universe conception of sets relies on this cumulative hierarchy of sets, and it is defined by transfinite induction on ordinals in the following way.
\vspace{0.5cm}

$V_0=\emptyset$,

$V_{\alpha+1}=\mathcal{P}(V_\alpha)$,

$V_\alpha=\bigcup_{\beta<\alpha} V_\beta$, if $\alpha$ is a limit ordinal.\footnote{$\mathcal{P}(V_\alpha)$ denotes the {\em power set} of $V_\alpha$, that is, the set of all subsets of $V_\alpha$.}

For all ordinals $\alpha$, $V=\bigcup_\alpha V_\alpha$.
\vspace{0.5cm}

All axioms of ZFC are true in $V$.\footnote{Of course the claim that $V$ is a model of ZFC cannot be proved from ZFC due to G\"{o}del's Second Incompleteness Theorem. Moreover, $V$ is proper class and such objects formally do not exist in ZFC set theory.} Any `set' that could be conceived is an element of $V$, so $V$ is {\em all} there is.\footnote{However, it is known that by forcing extensions one could add certain objects and make $V$ look `fatter'.} The {\em intended model} (or the {\em standard model}) of set theory is usually realized as a transitive class based on the iterative conception of sets, beginning with the empty set and taking the power set of sets formed in the previous stages.\footnote{A set $A$ is {\em transitive} if $x\in A$ implies $x\subseteq A$.} One way of thining down $V$ is by restricting the power set operation. G\"{o}del's {\em constructible universe} $L$ is obtained by restricting the power set operation to sets that are `definable' from the sets defined in the previous stages. So instead of considering {\em all} subsets in each successor stage, G\"{o}del considered subsets that are definable by a first-order formula. G\"{o}del also proved that $L$ is a model of set theory. In fact, $L$ is quite `small' in the sense that the existence of genuine large cardinals becomes inconsistent in $L$. Ideally, then, we want to find `enlargements' of $L$ so that genuine large cardinals can exist in these enlarged universes. This is roughly known as the {\em inner model programme}.\footnote{See Jensen \cite{Jensen} (1995), Steel \cite{SteelInnerModel} (2010), and Woodin \cite{WoodinUltimateL} (2017) for further reference on this topic.}

The realist philosophers of set theory who support monism in general argue for that, ultimately, every proposition of set theory must be either true or false in $V$. Essentially this is {\em epistemological realism} (or {\em realism in truth value}), which is the view that every mathematical statement is either true or false, independent of us, in an absolute Platonic universe.\footnote{See Shapiro \cite{Shapiro2000} (2000), \S 2.2, for more details. See also Wright \cite{Wright1986} (1986), p. 9.} G\"{o}del proved that if we add to ZFC the {\em Axiom of Constructibility}, that is, $V=L$ ---the statement that every set is definable---then we can prove CH. Cohen, on the other hand, constructed {\em forcing extensions} of a given model of set theory by adding `generic reals', and he showed that CH can be made false in these forcing extensions. In fact, we can also construct forcing extensions in which CH is true. Given a model of set theory, CH turns out to be true in some forcing extensions of the model, and false in other extensions. Many set theorists believe that both CH and $\neg\textrm{CH}$ universes are plausible models of set theory.\footnote{See for example Hamkins \cite{HamkinsMultiverse} (2012).} The existence of these kinds of universes leads to the pluralist conception in the set-theoretic arena.

G\"{o}del, as a Platonist, did not believe in the multiverse conception of set theory as he separated the situation we have with set theory from what we have with non-Euclidean geometries.\footnote{See G\"{o}del (1964), p. 483. Page references to the reprint version in Benacerraf and Putnam \cite{BenacerrafPutnam} (1983).} For G\"{o}del, there must be a fact of the matter as to whether CH is true or false. G\"{o}del justifies his claim relying on the apparent motivation that geometry concerns mathematical modelling of the physical space, whereas the interpretation of the set theoretical language is purely based on our mathematical intuition. 


Another reason that set theory cannot resemble geometry in this sense is the challenge to understand the set theoretical universe $V$ and how complex structure-wise it is compared to spatial structures in geometry. I will give a discussion about this more at the end of the paper.

\subsection{Multiverse Conception of Sets}\label{sec:multiverseconception}

The multiverse conception of sets is a major and new theme in the philosophy of mathematics. In fact, most of pluralist ideas about mathematics today revolve around the studies on the multiverse conception of sets.\footnote{The idea that there could be many set conceptions was, however, earlier anticipated by Mostowski and Kalmar in the mid 1960's. See Footnote 1 in Hamkins (2012) for further details.} Prominent supporters of this view in the mathematical community include Joel David Hamkins \cite{HamkinsMultiverse}, and Sy-David Friedman \cite{AntosFriedman}. For more methdological concerns, further studies can be seen in Davies \cite{DaviesPluralism}, who defends pluralism under Errett Bishop's \cite{Bishop} constructive approach of mathematical analysis. Hellman and Bell \cite{HellmanPluralism} argue for the necessity of pluralism at the foundational level in mathematics. Finally, Mark Balaguer \cite{BalaguerPluralism} argues for that both Platonism and anti-Platonism entail a version of pluralism. Readers may also find studies related to set-theoretic potentialism, which will be discussed later on.

The concept of set is usually understood, especially by non-set-theorists in general, as a single set concept belonging to a single set-theoretical universe, that is the cumulative hierarchy $V$. However, as we earlier said, results in set theory indicate that there are other universes of sets just as there are other geometries. The view that there are many universes of sets existing independently in the Platonic sense is known today as the {\em multiverse theory}. The idea that there is only one absolute universe of sets is called the {\em universe} conception.\footnote{See Hamkins (2012), p. 416, for his definitions of the universe and the multiverse conceptions.} Apart from the cumulative hierarchy $V$, other well known examples of set-theoretic domain include G\"{o}del's constructible universe $L$, the class of ordinal definable sets OD, the class of hereditarily ordinal definable sets HOD, and many others. All these classes admit different models for ZFC.

The notion of independence is generally regarded as a distracting issue for the adherents of the universe view.\footnote{I will point out later that independence is not as distracting as one might think in a hierarchical picture of multiverse of sets, as we can still preserve a sense of epistemic objectivity.} The multiverse conception of sets can bypass the infamous independence phenomenon more easily. According to Hamkins, the universe conception limits future insights that may arise from these possible mathematical realms.\footnote{We encourage the reader to refer to Hamkins (2012) for a detailed account of his position. For a criticism, we refer the reader to Koellner \cite{KoellnerHamkinsMultiverse} (2013).} 


How does the multiverse view settle the most fundamental questions of set theory which are known to be independent of ZFC? Let us take CH as an example. The problem with the universe view is that CH currently remains as an undecidable proposition. Nonetheless, for the realist, it has to be true or false in the reality of sets (such as $V$) that we yet to understand. But on the multiverse view CH is already settled, perhaps in a non-deterministic manner: CH is true in some universes, and `equally' false in others. Supporters of the universe view might suggest that CH would be solved using the following schema, which Hamkins refers to as the {\em dream solution}.\footnote{Hamkins, p. 430, 2012.}
\vspace{0.5cm}

Step (i). Produce an `obviously true' statement $\varphi$ of set theory.

Step (ii). Prove that $\varphi$ implies CH (or $\neg\text{CH}$).
\vspace{0.5cm}

Koellner claims that there has never been such a dream in this sense.\footnote{Koellner,  p. 24, 2013.} Even if there is so however, according to Hamkins, such a dream solution is impossible to achieve. This is due to the fact that mathematicians have had rich experiences in many various set theoretical universes.\footnote{But as we shall see, there are also arguments in favour of that some set theoretical universes are {\em richer} than others, or generally that some are more `preferred' over the others. In that case, the implications of the less preferred theories are also expected to be `less preferred' than that of the more preferred theories. If theory $T$ is preferred over theory $S$ and if $X$ proves CH while $Y$ proves the opposite, then it may be more reasonable to favour CH rather than the opposite.\label{footnote:dreamcounter}} Set theorists are now deeply informed about the status of CH and continuing their investigations both in the CH universes and the $\neg\text{CH}$ universes. For Hamkins, if we are to introduce an obviously true statement $\varphi$ which implies $\neg \text{CH}$, then we can no longer see $\varphi$ as an obviously true statement since doing so, would negate the experiences in the universes in which CH holds.\footnote{On hierarchical pluralism, either the experiences in the CH universes or the experiences in the $\neg\textrm{CH}$ universes weigh more. So from this viewpoint, negating the experiences in one universe is just decreasing its plausibility that was initially thought to be more plausible than our other experiences.} Similarly, if we are to introduce some $\varphi$ which implies CH, then the same argument would hold for the opposite.

Ternullo \cite{Ternullo2019} describes Hamkins' multiverse vision of allowing one to inspect set-theoretical universes from the point of view of the others as {\em Multiverse Perspectivism}. According to Ternullo, given $V$, the set-theorist is interested in knowing what universes one may have access to from $V$ and how we access them, and what $V$ looks like from the perspective of the other universes.\footnote{Ternullo, p. 56, 2019.}

The paper by Antos, Friedman, {\em et al} \cite{AntosFriedman} has also received attention in the multiverse community. The authors present their own multiverse theory. They define the universe and multiverse views similarly as in \cite{HamkinsMultiverse}. However, they introduce an additional criterion of differentiation between the conceptions. They call this criterion {\em commitment to realism}, which `measures how strongly each conception holds that the universe or the multiverse exist {\em objectively}'.\footnote{Antos, Friedman, {\em et al}, p. 3, 2015. At least commitment to some modest form of realism if not to a radical form. This is also congruent with the idea of having `degrees of intentionality' for models. However, the concern of the authors in \cite{AntosFriedman} is ontological rather than epistemological.} It was also mentioned in \cite{AntosFriedman} that the universe and the multiverse views would each split into two further categories according to whether one is a realist or a non-realist, yielding four positions in total: Universe view realism, universe view non-realism, multiverse view realism and multiverse view non-realism.\footnote{ibid, p. 3.} The authors decribe the universe view realism alternatively as G\"{o}del's realism. On the other hand, a universe view non-realist uses the set theoretical universe merely for practical purposes without actually believing in its independent existence. The two other positions are the multiverse view realism and the multiverse view non-realism. The adherents of the former view coincides with Balaguer's {\em full-blooded Platonism} (FBP), the view that any mathematical object that could exist actually does exist.\footnote{See Balaguer's \cite{Balaguer1998} and \cite{Balaguer2016} for details. When FBP is used in the form of ontological and epistemological pluralism, it constitutes the most radical version of the multiverse theory. I give a critique of such versions in the beginning of Section \ref{sec:mainsectionHP}.} Finally, a multiverse view non-realist does not believe in the existence of any universe, nevertheless she takes the multiverse phenomenon as a mathematical practice. Antos, Friedman, {\em et al.} also discuss a relatively new multiverse conception where they claim that Zermelo's \cite{Zermelo} account of the universe of sets can be further extended to what they call the {\em vertical multiverse} conception, in which they allow the cumulative hierarchy $V$ to be heightened while keeping the width fixed.\footnote{G\"{o}del's idea---that axioms should have a {\em maxizing property}---can be divided into two cases: ordinal maximality and power set maximality. That is, we want to maximize the size of the class of ordinals in the theory as well as the number of subsets of ordinals. Now $V$ is the `longest' and the `thickest' universe in this sense, as it contains {\em all} ordinals and {\em all subsets} of each ordinal. The height of the v-shaped cumulative hierarchy depicts the amount of `sizes of sets' in the universe, e.g., cardinals. Whereas the width shows how many, say, real numbers exist in each level of the cumulative hierarchy. From a universist/monist point of view, we cannot properly extend the height of $V$ by the forcing method since $V$ contains all ordinals. So Cohen's method instead expands the width of the cumulative hierarchy by adding generic reals. This is one reason why there are so many independent problems in set theory and it is not hard to define new statements that are independent from the axioms of ZFC. This is also the reason why set theorists desire to find some structure of the set-theoretical universe, that is, to be able to discover what there is within the structure, rather than arbitrarily defining what they are.} The authors also explain the distinction between what they call {\em actualism} and {\em potentialism}, two positions regarding the modifiability of $V$. Now, the idea of taking the cumulative hierarchy of sets as a `potential' universe rather than `actual' one was, in fact, discussed earlier by Linnebo \cite{Linnebo2013}.\footnote{See \cite{HamkinsLinnebo2019} for a study on the modal logic of set-theoretic potentialism. See also \cite{Berry2022}.} However, the foregoing study concerns our work to a greater extent. The authors in \cite{AntosFriedman} describe actualism as the belief that the cumulative hierarchy $V$ is a fixed object and that it is impossible to `stretch' $V$ via model theoretic constructions.\footnote{Antos, Friedman, {\em et al.}, p. 19, 2015.} Any forcing extension of $V$, according to actualism, is in fact a model in $V$. The opposite extreme is called {\em potentialism} which they define in the following manner. 

\begin{changemargin}{1cm}{1cm}
{\em A potentialist, on the other hand, sees $V$ as an indefinite object, which can never be thought of as a `fixed' entity. The potentialist may well believe that there are some fixed features of $V$, but she believes that these are not sufficient to fully make sense of an `unmodifiable' $V$: the potentialist believes that $V$ is indeed `modifiable' in some sense.\footnote{ibid, p. 19.}}
\end{changemargin}

A major section in \cite{AntosFriedman} is reserved for the explanation of how the vertical multiverse view is used in their {\em hyperuniverse programme}, which I will discuss shortly in more detail in relation to the hierarchical setting of the multiverse.

\section{Hierarchical Pluralism}
\label{sec:mainsectionHP}
I shall first begin with a criticism of radical versions of mathematical pluralism or radical multiversism by asking the following provoking question: Does pluralism liberate the mathematical ontology or blur mathematical truth? It appears to me that there is trade-off, and that we ideally try to find a balance between the following two extremes: the richness of mathematical ontology and the objectivity of mathematical epistemology. I shall call this dichotomy the {\em pluralist's dilemma}. According to this dichotomy, it is not possible to fully satisfy either side without compromising the other. In fact, one of the main motivations in establishing the view presented here is a natural consequence of the desire to find an equilibrium between the two extremes. On the one hand, we cannot deny the plurality of set-theoretical universes. This is not just due to the experiences that set theorists had and many fruitful results they obtained in these different universes, but it is also that some of them can be regarded as {\em intrinsically} plausible, and that these models cannot be easily singled out to a unique absolute model. On the other hand, radical forms of pluralism, to a certain extent, are rather too liberal about the epistemology of mathematics, which overcomplicate settling the truth value of statements.\footnote{For more ontological concerns, see Barton \cite{Barton} (2016).} The main worry is that, as tempting as it may sound to liberate mathematical truth, we depart from objectivity in determining even intrinsically justifiable statements of mathematics. Radical pluralism must reject mathematical intuition. According to this type of pluralism, there are in fact many mathematical intuitions possibly incongruent with each other.\footnote{A radical realist may think that liberating mathematical truth and ontology leads to paraconsistent mathematics with no basis for a credible assessment of statements.} If two contradictory statements can be satisfied in two different universes respectively, then what determines a sensical mathematical intuition? For a pluralist of this kind, there is in fact no common mathematical intuition. A radical pluralist has no or, at best, a little motivation in finding what is intrinsically appealing, nor can she rely on any legitimate reason for her objection to `counter-intuitive' statements being called `true' statements. This is due to that a radical pluralist---{\em Hamkinsian pluralism}, one might say---is positioned on one extreme of the pluralist's dilemma: richness of mathematical ontology. But this is reinforced by adopting the philosophy that we should have as many variety of set theories as possible rather than having intrinsically justified statements so as to favour the other extreme of the pluralist's dilemma, namely the objectivity of mathematical truth. On what grounds do radical pluralists decide to eliminate counter-intuitive statements? For a radical pluralist, there is a very good reason not to believe in the Axiom of Extensionality equally as to believe in the existence of rank-to-rank embeddings, and even ones which contradict ZFC, e.g., Reinhardt cardinals. Furthermore, both beliefs can be considered as {\em equally} plausible by the radical pluralist. It is generally held that the structure of natural numbers {\em is} the intended model of arithmetic.\footnote{It may be said that this is begging the question for monism. But it seems to me that the natural number structure supports monism only for arithmetic, not for a more general concept like `set', whatever its structure is. The reason there is a consensus on intended arithmetic but not on the concept of `set' is that I think arithmetic has more structural and informational content than sets have. So the idea that the structure of naturals is the intended interpretation of arithmetic does not seem to favour monism regarding sets, but it may, at best, promote a monist view of a specific segment of the cumulative hierarchy whose elements are finite sets.} This is of course due to that the structure of naturals is considered to be intrinsically plausible by so many mathematicians. But how do we decide if another structure for arithmetic is not as plausible as the natural numbers? Even if the case for arithmetic is more obvious, the latter question makes things more challenging in set theory as the concept of set is not quite as `uniquely' interpreted as natural number arithmetic. Despite the fact that there are many set conceptions, even for a modest realist it seems inaccurate to put two different conceptions or theories of sets on the same level of intention. 

The liberation that pluralism offers the mathematician, when taken in the extreme sense, comes from the expectation that pluralism could as well be applied to the methodology. In \cite{PriestMathPluralism}, Graham Priest lists a number of objections raised by anti-pluralists and replies to them. One particular objection concerns the universality of the inference rules in proving mathematical statements. Priest says:

\begin{changemargin}{1cm}{1cm}
{\em $[...]$ when someone breaks the rules of a mathematical game, and they are not simply making a mistake, they are, {\em ipso facto}, no longer playing that game. The new set of rules constitutes a new game.\footnote{Priest, p. 9, 2013.}}
\end{changemargin}

The set of rules of inference is one thing that any non-formalist mathematician would want to be certain about its soundness. This type of methodological pluralism, i.e., pluralism applied on the derivational rules of the mathematical game, leads one to have non-standard interpretations of the logical language. Since some primitive rules such as the laws of thought or Peano axioms, etc., are more universally accepted than otherwise, there is a good reason these axioms or rules deserve to be called `standard'. It seems to me that the reason why they are universally accepted is not because of any pragmatic convention or that they have been taught to us that way, but it is because of how primitive and finitary facts about quantities, collections, space, and so on, are conceived by us via empirical observations.\footnote{For example, it appears to me that the reason we have `$x+0=x$' as an axiom of arithmetic is mostly due to the fact that we cannot find an empirical counter-example to the sensory observation that adding nothing to a collection yields the collection itself. So `$x+0=x$' is just an abstraction that generalises this empirical observation. I think that the more abstract the objects get, the more we tend to rely on the usefulness of these objects to the mathematical discipline for their existence. For example, infinite dimensional vector spaces or infinite sets may not correspond to any empirical observation, yet we usually accept and define them mostly for pragmatic reasons. Now despite that infinite dimensional vector spaces may not correspond to any our empirical observations, the smaller components that constitute an infinite dimensional vector space eventually are constructed from abstracting our sense data. For instance, a hypercube itself will not correspond to an abstraction of our empirical experiences about the physical world, but its components (cubes, lines, intersecting two points with a line) will.}

Due to concerns given above, I do not endorse radical versions of the multiverse theory for the fact that radical pluralism is inevitably inconsistent with even modest versions of realism, especially in the epistemological (truth-value) sense. The idea is then to adopt various sorts of criteria for comparing any two competing theories for sets so as to decide which one is more plausible and is a better approximation to the `full conception' of sets and then look at the statements in the most plausible theory.\footnote{I am using the term {\em full conception} to describe the `true theory of sets', i.e., true statements in the Platonic concept of `set'. We know that $V$ is all there is. The problem though is to actually classify $V$ in a more meaningful way. There are theories that suggest for example that $V=L$, or that $V=\textrm{HOD}$, and more recently Woodin's $V=\textrm{Ultimate-}L$. For the latter, see Woodin \cite{WoodinUltimateL} (2017).} The existence of full conception is supported by the preservation of {\em foundationalism}, as Maddy points out that `if set theorists were not motivated by a maxim of this sort, there would be no pressure to settle CH'.\footnote{Maddy, p. 209, 1997.} Now one might argue that the foundational theory does not need to comply with the full conception in reality.\footnote{At this point, I think I should also draw a line between {\em foundationalism} and the idea of {\em reductionism}. Reductionism refers to the idea that mathematical objects can be reduced to objects of {\em one of possibly many} simplified theories. Foundationalism, on the other hand, is the view which argues for the existence of a {\em unique} foundational theory for all (or most) of mathematics. It is unclear whether one should consider the uniqueness requirement as a part of the definition of foundational theory, but I think uniqueness has to be added so that the foundational theory serves its purpose correctly. I think it is safe to say that reductionism does not require this uniqueness condition. Hence, we may think of foundationalism as a strict form of reductionism.} But it must be congruent with how we understand the full conception and the way it appears to us. 

Now, although there may be a pressure to settle CH for the foundationalist, in the geometric scene there is no pressure to settle the parallel postulate. I think the reason for this is the following. It seems quite probable that Euclid's fifth axiom must have been posited in the Euclidean conception of space and in the way we experience physics, which was perhaps the only conception of geometric space at the time of Euclid. Yet it is not clear in what kind of structure CH was meant to be formulated by Cantor. The fact that non-Euclidean and Euclidean distinction is far clearer than the distinction between the CH and non-CH universes of sets led the majority of mathematicians to accept the fact that the parallel postulate is associated with a unique space context, that is, the Euclidean geometry. The two space conceptions are obviously different. However, the blurriness between CH and non-CH universes and that set theory has no natural structure inevitably leaves some room for research to find out whether CH is true in some `natural' conception of set that monists wish to describe. Now it is an interesting question for its own sake why geometry and set theory are so different in this sense. There is (or was) a programme to find a structure for sets, but there has been no similar motivation to find a natural conception of `space'. I think the desire to settle CH goes hand in hand with the quest of finding a structure for set theory. I do not think there is such an intention for geometry, as Euclid's {\em Elements} is a clear example to a material axiomatic system, where Euclid had a definite space conception in mind prior to the axioms. Thus, mathematicians leave out the parallel postulate and instead examine statements of set-theoretical nature.

Pluralist philosophy of mathematics can either be interpreted based on realism or nominalism. The pluralism I advocate for---that is hierarchical pluralism---is a view compatible with a form of realism in truth value, which also allows a multiverse interpretation of sets. Not every pluralist philosophy needs to be based on realism, however. For example, Balaguer's pluralist account is known to be a fictionalist view.\footnote{See Balaguer (2009), (2017).} Hence, just as there is no reasonable argument in favour or against realism/fictionalism---as Balaguer puts it---there is no reasonable argument in favour or against the existence of the full conception of sets. I will then take the realist route and adopt the idea that the full conception of set theory exists independently. I shall also endorse epistemological realism (realism in truth-value), that every mathematical statement---about set theory in particular---is true or false with respect to the full conception. I argue that the preservation of epistemological realism while maintaining the multiverse of sets leads to a more hierarchical picture of set-theoretical universes.

At this stage we can ask if one could compare two competing theories based on their plausibility or intendedness, as if having a `degree of intentionality' for each of these theories.\footnote{For example, in the case of natural number arithmetic, it is generally thought that Peano arithmetic is more intended than non-standard theories of arithmetic. So it is natural to put the intendedness of Peano arithmetic `above' that of non-standard theories of arithmetic, and so the consequences of Peano arithmetic will be more plausible than that of the given non-standard theory. This is not to claim that one consequences are true and the others are false, but it is to say that true statements of a more plausible theory have a higher sense of preferability than the ones in non-standard theories.} I will shortly explain what I mean by this. An immediate objection that could be raised right now is that this kind of epistemic assessment and comparison cannot be made at all. But I think this is wrong and it is not something that causes a problem. As a matter of fact, the idea that some set theories could be more plausible than others was even left open by Hamkins, who is known to support more liberal versions of the multiverse theory, at least in the sense that we are not obliged to assume that every universe in the multiverse should have equal {\em status}.\footnote{See Magidor \cite{Magidor2012} (2012) for a discussion that some set theories can be selected as more preferable than others.} Hamkins says that `we may prefer some of the universes in the multiverse to others, and there is no obligation to consider them all as somehow equal'.\footnote{Hamkins, p. 417, 2012.}

I would specifically like to emphasise on the consequence of the last sentence, quoted by Hamkins. It is worth noting that I do {\em not} advocate for any kind of {\em ontological} privilege or status of particular universes over the others. No universe comes ontologically prior to another. The {\em status} I refer to is rather epistemological, which may be granted either based on extrinsic or intrinsic justifications. So when I use the term preferable theories or universes, it should be understood that the concern here is epistemological rather than ontological.

Ever since G\"{o}del's programme, studying the methods of settling independent problems continued to attract many mathematical logicians and philosophers of mathematics. In the contemporary scene, for example, John Steel on the multiversist camp developed a maxim for determining mathematical truth. He suggested that `the key methodological maxim that epistemology can contribute to the search for a stronger foundation for mathematics is: {\em maximize interpretative power}'.\footnote{Steel, p. 154, 2014.} Steel's main idea is to maximize the interpretive power of the language of set theory by extending it to the {\em multiverse language} so as to `avoid' asking questions like CH. Among the multiverse studies, Shelah's {\em Logical Dreams} \cite{Shelah} is particularly interesting for that he proposes certain criteria for choosing axiom candidates. He argues that one could adopt a measure-like metric in determining which axioms are `true' axioms of set theory.\footnote{Throughout the paper I will use the term `metric' to refer to any criteria for ordering theories so as to define their relative plausibility distance in that sense.} He writes:

\begin{changemargin}{1cm}{1cm}
{\em First and most important, it must have many consequences, making it have a rich, deep, beautiful theory. Second, it is preferable that it is reasonable and ``has positive measure". Third, it is preferred that it leads to no contradiction $[...]$.\footnote{Shelah, p. 214, 2003.\label{footnote:Shelah}}}
\end{changemargin}

The {\em hyperuniverse programme}, originally proposed by Arrigoni and Friedman \cite{ArrigoniFriedman}, has similar purposes, and it is primarily associated with an on-going project of searching new axioms for set theory. Now, I should emphasise that it is {\em not} my intention here to give specific criteria for selecting preferred theories, as these have been already studied by Magidor \cite{Magidor2012}, Shelah \cite{Shelah}, and in the hyperuniverse programme; it is rather to argue about the consequences of introducing the concept of `preferred' theories or universes and discuss how this classification yields a hierarchical treatment of mathematical truth in a hierarchical picture of the multiverse if we want to preserve a considerable amount of epistemic objectivity about mathematical truth, i.e., to find an equilibrium in pluralist's dilemma.\footnote{For a philosophical discussion about preferable models, see Button and Walsh (2018), pp. 39-44.} However it is also worth looking at how the hyperuniverse programme handles the problem of selecting set-theoretical universes. The hyperuniverse programme is one where the authors give a set of general requirements for preferred theories by what they call the {\em omniscience} and {\em maximality} conditions. One point where hierarchical pluralism differs from the hyperuniverse programme is reference to Platonism. Platonism is not invoked in the hyperuniverse programme, whereas the existence of full conception of sets is by itself a Platonistic reference, though this reference does not go beyond talking about the true statements that hold in $V$. This will not put any limitation on the treatment of hierarchical pluralism, nor will it ignore the multiverse of sets existing independently. Another apparent difference is that hierarchical pluralism is actually a realist-based view compatible with the idea of many-universes, whereas the hyperuniverse programme seems to eliminate foundationalism for the fact that it avoids to make any Platonic reference and so any selected theory in that sense has no relation to the full conception of sets. So I intend to introduce hierarchical pluralism as a result of selecting preferred universes over the others and discuss how truth can be determined in such a dynamic picture. The preferred theories, regardless of how they are selected, must also admit epistemic priority over the non-preferred theories. This is where the hierarchical treatment of the multiverse comes into consideration. If two distinct set theories are comparable by any of the approaches introduced in the aformentioned studies, then one theory must be more plausible than the other. In other words, hierarchical pluralism is a natural product of deeming some theories more plausible than others.

So what can be the justification for having an ordered hierarchy of theories of sets? There is no unique justification for the existence of an ordering; in fact, there may be many. An ordering surely can exist. For example, one natural way of determining the status of theories could be by looking at how they are placed in the consistency strength hierarchy. Suppose that $T$ and $U$ are two theories extending ZFC. Then $T\leq_{con} U$ iff ZFC proves $\text{Con}(U)\rightarrow \text{Con}(T)$.\footnote{$\text{Con}(T)$ stands for ``$T$ is consistent", i.e., no contradiction can be derived from $T$.} If $T\leq_{con} U$ and $U\leq_{con} T$, then we write $T\equiv_{con} U$, and say that $T$ and $U$ have the same consistency strength or they are {\em equiconsistent}. There is a general consensus that consistency strengths of `natural' extensions of ZFC are linearly ordered.\footnote{For supporting views, see Caicedo \cite{Caicedo2011} (2011), Koellner \cite{Koellner2011} (2011), Simpson \cite{Simpson2009} (2009), Steel \cite{Steel2014} (2014). For an opposing view, see Hamkins \cite{Hamkins2022} (2022).} Steel says:

\begin{changemargin}{1cm}{1cm}
{\em If $T$ is a natural extension of ZFC, then there is an extension $H$ axiomatized by large cardinal hypothesis such that $T\equiv_{con} H$. Moreover, $\leq_{con}$ is a pre-well order of natural extensions of ZFC. In particular, if $T$ and $U$ are natural extensions of ZFC, then either $T\leq_{con} U$ or $U\leq_{con} T$.\footnote{Steel, p. 157, 2014.}}
\end{changemargin}

We know that due to G\"{o}del, there is no largest consistency strength (except of inconsistent theories). One can increase the consistency strength of a theory $T$ by adding to it the sentence that $T$ is consistent, or the statement  ``There exists a model of $T$". There is a good evidence that the plausibility of extensions of ZFC with greater consistency strength is in fact lower. Regarding this, Shelah writes:

\begin{changemargin}{1cm}{1cm}
{\em If you go higher, up in the large cardinal hierarchy, the justification for their existence is decreasing $[...]$ $($so lower consistency strength is better$)$.\footnote{Shelah, p. 214, 2003.}}
\end{changemargin}

Applying the maximize property in the G\"{o}delian sense to select preferred universes relies on the naturalistic approach. However, this dismisses the fact that the consistency strength of various large cardinal axioms are higher than others. So there is no {\em a priori} reason to consider different theories all equal in terms of their consistency strength. Thus from this viewpoint, the existence of rank-to-rank embeddings should be less plausible than the existence of inaccessible cardinals. So there is at least one example where `natural' extensions of ZFC are thought to be linearly ordered, and that is with respect to their consistency strength. If $\text{ZFC}+\varphi$ and $\text{ZFC}+\psi$ are two natural extensions of ZFC, then either $\text{Con(ZFC}+\varphi)\rightarrow \text{Con(ZFC}+\psi)$ or $\text{Con(ZFC}+\psi)\rightarrow \text{Con(ZFC}+\varphi)$. If one adopts consistency strength as her plausibility metric, then either $\text{ZFC}+\varphi$ or $\text{ZFC}+\psi$ is more plausible than the other.\footnote{It may also be thought that the consistency strength hierarchy is well-ordered. Hamkins (2022), however, claims that it is ill-founded.}

We start by considering all universes of sets. For example, in our case, the multiverse may contain ill-founded universes. Desideratum 1 in Arrigoni and Friedman \cite{ArrigoniFriedman} assumes that well-founded universes have a higher plausibility than that of ill-founded universes. We also generalise this comparison to be made for universes in which the axioms of ZFC hold. So some well-founded universes are more plausible than other well-founded universes. Independent of the method one may follow for selecting preferred theories or universes, if not all set theories have equal status in the multiverse, then some must {\em a priori} be more plausible than the others.\footnote{This at least leads to a more objective evaluation of theories and the truth value of set-theoretical statements in comparison to more radical versions of the multiverse theory.} By {\em status}, we mean `intrinsic plausibility' or `fruitfulness', depending on whether we endorse, respectively, an intrinsic position or a more naturalistic view. This is why I will encapsulate both under the term `intendedness'. Our intentions about using a concept in mathematics could be based either on intrinsic or extrinsic grounds, or mixed. But I think it is safe to say that our pragmatic intentions about a concept must be congruent with properties that are intuitive to us.\footnote{cf. Shelah's second requirement in our Footnote \ref{footnote:Shelah}.} These two properties, thus, are not entirely disjoint from each other. Here is a short argument why. As much as we take extrinsic justifications in determining the plausibility of a statement, we must also consider intrinsic justifications into account. Of course we also require some practical implication, as solely relying on intrinsic justifications does not `extend' our mathematical knowledge. In fact, logical tautologies are of this kind. All tautologies, by virtue of themselves, are intuitively true and so they are intrinsically justifiable. But they do not extend our mathematical knowledge, and they certainly do not entail too many statements. More specifically, tautologies only entail other tautologies. On the other extreme we have extrinsic properties, in the sense of abundance of logical implications. So another option when determining plausibility is to take just extrinsic justifications into account without looking at the intrinsic properties. Since anything follows from inconsistency, a contradiction is technically the most `logically abundant' statement one can have. Of course, we do not want our theories to be inconsistent, so we must necessarily take intrinsic justifications into account as well. This short argument is to show that there should be a balance between our intrinsic and extrinsic justifications. We should neither solely rely on the intrinsic justifications nor on extrinsic justifications of a statement to see if it is a natural axiom candidate. If any priority had to be given between the two types of justifications, my personal view would be however to endorse the motto ``intuitiveness before practicality".

So the idea is to have a plausibility ordering between the theories (or models) for a given concept. Whether our intentions are  grounded on intrinsic or extrinsic features is not relevant to our discussion though. The important point is that some theories---particularly theories about sets in this case---are expected to be more intended compared to others. If we look back to Hamkins' {\em dream solution} argument, it is possible to see that the argument is based on the assumption that both $\text{CH}/\neg \text{CH}$ universes have the same `degree of intentionality'. This is what I call the {\em equal-status interpretation} of the multiverse.\footnote{According to this interpretation, any two extension of ZFC are {\em equally} plausible. For example, $L$ is a model of set theory in the same plausible sense that HOD is. Whether we consider models or theories, the equal-status interpretation deems any two theories or any two models equally plausible. The problem with this interpretation is the following. Even if we have two extensions of ZFC, say $\text{ZFC}^+$ and $\text{ZFC}^*$, from which CH  and $\neg\text{CH}$ are provable respectively, there is no way to settle independent problems like CH with this interpretation for Hamkins' worry in the {\em Dream Solution} that choosing one over the other would negate the experiences we had in the other theory. On hierarchical pluralism one should be able to settle CH and other independent questions {\em within} the `most plausible theory' in hand.} But as Hamkins says, `there is no obligation to consider them all as somehow equal'. So an independent position in the multiverse view can be that some of these theories may have a `higher' priority. 
 Regardless of what our preference metrics are, it is possible then to have an ordering between them with respect to the chosen metric.

Again, classifying the universes with respect to their degree of intentionality can be done in two ways: Either by a naturalistic approach, in which case by looking at the fruitfulness of the universes (that is, relying on the extrinsic justifications), or by purely based on the intrinsic justifications. I do not advocate for a particular approach in this paper, nor do I propose a specific method on how to determine if a statement is intrinsically plausible. However, since any modest version of realism entails that there exists a full conception of sets, any two theories about sets are {\em in principle} comparable in terms of intrinsic plausibility with respect to the full conception. In any case of conflict between two competing obvious statements, it can still be decided in principle which one of two competing statements is more obvious than the other. Following Hamkins' own example, if someone proposes an obviously true statement $S$ implying that Manhattan does not exist, then we decide whether the existence of Manhattan or the truth of $S$ is more obvious. If the existence of Manhattan is more obvious than $S$, then it must be that our base theory plus the existence of Manhattan is a better approximation to the full conception of the theory in consideration. In case of otherwise, our base theory plus the truth of $S$ is a better approximation. 

We said earlier that the equal-status interpretation of the multiverse is epistemologically problematic as it leads to radical multiversism. Claiming that the `degrees of intentionality' of different theories are on a par with each other does not help us settle independent statements. It is by some sort of classification of these theories or models we get to pick the best possible one (in this case, the one with the `maximal' degree of intentionality---as I shall define it shortly), and regard that as the best possible reference. It is expected that this theory may later get replaced in the intentionality hierarchy by a better one. Our best possible set theory naturally forms the foundational theory. So one feature that can be seen in hierarchical pluralism is Mich\`{e}le Friend's concept of `growing foundational theory'. In \cite{FriendPluralism}, Friend gives four anti-foundational arguments one of which is based on the idea that the foundational theory is a continually {\em growing} theory. She says: 

\begin{changemargin}{1cm}{1cm}
{\em $[...]$ the foundationalist begins with the technical result that most of mathematics can be reduced to The Foundation. This is a twofold mis-description. First, the reduction is sometimes too contrived, and therefore, not successful. Second, any proposed foundation is only that: a foundation. That is, we can add more to the foundation. Whatever the founding theory is,`it' grows. As it grows we understand the founding parts in a new light. So it is not a fixed foundation.\footnote{Friend, p. 60, 2014.}}
\end{changemargin}

Friend's {\em growing foundation} seems to be quite compatible with the naturalistic idea that the best possible theory is always improved, as if we approximate to the full conception. An extension of ZFC with the `maximal' degree of intentionality may have a growing feature for that reason.

Next is to discuss how truth can be determined in a hierarchical picture of the multiverse of sets. It was claimed by Friend that `foundational pluralism' is an unstable position.\footnote{ibid, pp. 24-25. For Friend's definition of foundationalism, see ibid, p. 8.} However, this might only be the case if we were to abandon naturalism. This is because of the fact that naturalistic approach to mathematics, in finding a structure for set theory, leads to a dynamic picture. For example, suppose that we accept a large cardinal axiom that we thought initially was fruitful. It may later turn out that we choose to withdraw this axiom in favour of another conflicting axiom. So then our foundational theory about sets is subject to change and improve itself. Friend's claim also seems to depend on whether we endorse the equal status interpretation or not. If all universes are plausible on a par with each other, then this type of radical pluralism may not be compatible with foundationalism. However, if the plausibility of some theories are higher than others, then it is natural to consider our best theory as our foundational theory. Now it still holds that truth only makes sense in a structure.\footnote{In fact, on hierarchical pluralism, one can go one step further and argue that truth makes more sense in our best possible theory about sets.} Friend calls this ``On truth in a theory" and she presents this as an anti-foundational property. She writes:

\begin{changemargin}{1cm}{1cm}
{\em Together, the structuralist and the pluralist do not think that there are absolute truths in mathematics of the form: ``\ $2+8= 10$". Instead, what is true is: ``In Peano Arithmetic, $2+8 = 10$".\footnote{Friend, p. 67, 2014.}} 
\end{changemargin}

More specifically, it seems that {\em if} Peano Arithmetic is our best possible theory for describing the natural number structure, then $2+8=10$ holds in Peano Arithmetic. Also, I think mathematicians would be willing to accept Friend's view as a kind of {\em local realism}.\footnote{I coined the term {\em local realism} to refer to a two-tiered temporal ontology which can be defined as follows. Mathematical statements are independently true or false in structures (and nowhere else) {\em once} the structures are defined. For example, on local realism one cannot tell whether or not natural numbers have a maximal element, as there is no fact of the matter whether or not there is a greatest natural number. But once the natural number structure is defined, the statement becomes `independently' true or false in that structure. The truth value of statements are not---in the temporal sense---independent of mathematical practice but it becomes independent once the background ontology is constructed. Local realism in this sense may be associated with Maddy's {\em thin realism} with temporal elements involved. See Maddy \cite{Maddy2011DefendingAxioms} (2011), Part III.} 

Before we elaborate more on the epistemological features of hierarchical pluralism though, we shall also look at Balaguer's intention-based Platonism (IBP) and his claim that IBP entails a type of relativism in mathematics.\footnote{Most of Balaguer's ideas on IBP can be found in his \cite{Balaguer2009}. In \cite{BalaguerPluralism}, Balaguer argues that IBP is compatible with pluralism.} Balaguer claims that both Platonism and anti-Platonism entail a certain version of pluralism. He then argues that the entailed version of mathematical pluralism is actually true. Considering the fact that both CH and $\neg\text{CH}$ are consistent with ZFC, Balaguer asks what Platonists should say about the truth value of CH. He mentions two types of Platonisms, namely {\em Silly Platonism} and {\em Better Platonism}. According to Balaguer, mathematical truth in Silly Platonism is determined as follows:

\begin{changemargin}{1cm}{1cm}
{\em A mathematical sentence or theory is true just in case it accurately characterizes some collection of mathematical objects. Thus, since the mathematical realm is plenitudinous, it follows that all consistent mathematical theories are true. And so it follows that CH and $\neg\text{CH}$ are both true, because CH is true of some parts of the mathematical realm, and $\neg\text{CH}$ is true of others.\footnote{Balaguer, p. 383, 2017.}}
\end{changemargin}

It is not this kind of platonism Balaguer relates with pluralism though. Silly Platonism is a very loose version of realism in which the theory could be realized in non-standard models, where the sentences satisfied in those non-standard models are still said to be true on a par with the sentences satisfied in the standard models. We see that the role of intended (or standard) models becomes more critical in his Better Platonism for which he describes mathematical truth in the following manner:

\begin{changemargin}{1cm}{1cm}
{\em There is a difference between being true in some particular structure and being true simpliciter. To be true simpliciter, a pure mathematical sentence needs to be true in the intended structure, or the intended part of the mathematical realm---i.e., the part of the mathematical realm that we have in mind in the given branch of mathematics.\footnote{ibid, p. 383.}}
\end{changemargin}

According to this view of Platonism, an arithmetical sentence is true iff it is true in the `standard model' of arithmetic. Balaguer claims that the intended structure may not be unique in every branch of mathematics, in particular set theory. For this, he gives an example of two structures $H_1$ and $H_2$ in which, respectively, $\text{ZFC}+\text{CH}$ and $\text{ZFC}+\neg\text{CH}$ hold. Balaguer also assumes that both $H_1$ and $H_2$ could be intended models of set theory as they could be both {\em fully} consistent with our full conception of sets.\footnote{See Balaguer (2009), pp. 142-144, for more details on his definition of the `full conception'.}

To extend Balaguer's view, one may further break down the simpliciter truth for preferred theories between themselves and argue something along the following lines: There is a difference between being true simpliciter in more preferred theory of sets, and being true simpliciter in a less preferred theory of sets. Both are true simpliciter, yet one has a higher sense of universality. So regarding Balaguer's account on intended models, an intriguing question that can be asked is whether $H_1$ and $H_2$ can really be {\em equally} or {\em fully} consistent with our full conception of set theory, and hence whether both can be taken as intended models {\em on a par} with each other, considering that they satisfy mutually contradicting statements.\footnote{This is one aspect of hierarchical pluralism that contrasts with Balaguer's IBP, as I will mention shortly, that if there are two competing theories, they should not be `equally' consistent with the full conception. But since set theorists have experiences in both universes, the multiverse view requires us to accept both universes while keeping in mind on hierarchical pluralism that one has a higher sense of intentionality than the other.} If these structures satisfy mutually contradictory statements, then to what extent a structure is said to capture our full conception of sets? Unless our full conception is contradictory, it seems that one of them must be `less' standard than the other. Intentionality of models is a varying concept. Some structures may serve as a `better' intended model than others, and not all structures---up to isomorphism---have the same level of intentionality. This is implicit in Shelah's {\em Logical Dreams}, apparent in Magidor's work, the hyperuniverse programme, and it is even left open by Hamkins himself who is known to support more liberal versions of the multiverse theory. Assuming an {\em equivalency in intentionality} (or the {\em equal-status interpretation} of the multiverse), otherwise, would make the line even between standard and non-standard models rather blurred. So on a hierarchical multiverse setting, not only can a distinction be made between standard and non-standard theories, we should also be able to order standard theories among themselves.

The question that what determines a structure to be counted as an intended model was also considered by Balaguer himself, and in fact he gives a criterion for this.\footnote{See Balaguer (2009), p. 144, for details.} I think the question is more complicated than it looks. Instead, let us consider the negation of the same question: What determines a structure {\em not} to be counted as an intended model? Apart from assigning an unconventional interpretation for objects and the language, the idea we discussed in the previous paragraph may give an answer to the negated version of the question. In fact, this does not settle which structure actually {\em is} or {\em is not} an intended model. It settles a rather weaker question. We can say, if $H_1$ and $H_2$ are two structures such that $H_1\models \varphi$ and $H_2\models\neg\varphi$ for some statement $\varphi$, then at least one of them cannot be an intended model `to the degree' that the other can. This is not to say that neither is an intended model. It just says that both structures cannot be simultaneously counted as intended (standard) models on the same level. The reason is that if $\varphi$ is a statement which is modelled by some standard structure of a theory and if $\neg\varphi$ is modelled by another standard structure, then there should be no reason to believe in existence of a consistent full conception, in which case this would completely abandon realism. For mathematical pluralism and epistemological realism to fit together, it appears that every model should have a `degree of intentionality' with respect to the full conception, so that some models are meant to be `more intended' or `more preferred' than the others. The view I express here is not a Hamkinsian form of pluralism---which relied on the equal-status interpretation of the multiverse---but it rather admits a relational type of pluralism which we may call {\em hierarchical pluralism}, where theories or models about a concept are `ordered' with respect to their intendedness. A hierarchical pluralist is an epistemological realist, i.e., realist in truth value. She might believe in the multitude of set theories, but she also believes that these theories are ordered with respect to their degree of intentionality---by fixing any of the metrics introduced by Shelah, Magidor, or the conditions provided in the hyperuniverse programme. The purpose of these metrics of course is to measure to what extent the theory conforms to the full conception of sets. In fact, for any tenable version of realism, there must exist such an ordering of the theories, as this is the foremost effect of epistemological realism on the multiverse theory. A hierarchical pluralist believes that sentences that are provable from theories with higher degree of intentionality have a higher sense of plausibility, preferability, or rather higher sense of truth, perhaps by assigning it a degree of truth in accord with its degree of intentionality.

Hierarchical pluralism compromises on elements of both radical forms of realism and of pluralism. It does not use a direct reference to absolute truth when defining truth. On the contrary, in the desiderata given below, truth will be defined dynamically in such a way to conform our best possible known theory as opposed to a fixed Platonic theory of `true set theory' or `true arithmetic'. Also, hierarchical pluralism is not a radical form of pluralism as one might expect. The idea of having degrees of intentionality for theories is what separates and classifies distinct theories in accordance with their plausibility. Radical multiversism has no concern to put theories of sets in order of preference, as they all can be taken as equally plausible on such a view.

I now give four desiderata of hierarchical pluralism which will shape the view better. It is best to begin with stating the existence of a metric to measure the plausibility of theories.\footnote{A conception can be described either by a model or by a theory. A theory may have more than one model. Given a model $M$, the set of true sentences of $M$ gives us the {\em theory of $M$}. So in our case, the ``true theory of sets" or the ``full conception of sets" is really just the set of sentences that are true in $V$. We may use {\em theories} and {\em models} interchangably here as the distinction between the two does not matter when speaking about plausibility.} This could be determined by one of the available methods, e.g., Shelah \cite{Shelah}, Arrigoni and Friedman \cite{ArrigoniFriedman}, or looking at the consistency strength. This then defines a degree of intentionality (or degree of plausibility) for each theory.
\vspace{0.5cm}

\noindent{\em 1. Degree existence.} Every model up to isomorphism or theory of a mathematical concept (e.g., set, group, arithmetic, field) has a {\em degree of intentionality}, which is a measure for its plausibility with respect to the full conception of the theory in consideration. 
\vspace{0.5cm}

As discussed earlier, we want to choose a metric that defines what it means for a theory to be more plausible than the other. The evaluation should assign every theory an associated degree of how plausible the theory is. The higher its plausibility, the more likely it becomes our best possible theory and the better approximation it is to the full conception. Determination of the plausibility of a theory is twofold. It could either be relied mostly on extrinsic or intrinsic justifications. If we lean more towards the intrinsic side, theories with higher consistency strength will be less preferable since they are closer to inconsistency. In case of otherwise, if extrinsic justifications are taken as our criteria for plausibility, theories with higher consistency strength will be more preferable as they give a more fruitful set theory in the sense of abundance of ontology. Either way one should be able to assign a plausibility degree for each theory based on the the decided metric.
\vspace{0.5cm}

\noindent{\em 2. Linearity}. Let $T$ be a theory for a concept with a degree of intentionality ${\bf a}$, and let $T'$ be a theory for the same concept with a degree of intentionality ${\bf a'}$. Then ${\bf a\leq a'}$ or ${\bf a'\leq a}$. If there exists a sentence $\varphi$ such that $T\vdash\varphi$ and $T'\vdash\neg\varphi$, then ${\bf a \neq a'}$. A theory with a higher (greater) degree of intentionality is accepted to be a better approximation to the full conception. Intuitively, the degree of intentionality ${\bf b}$ of a theory $S$ is higher than the degree of intentionality ${\bf a}$ of the theory $T$ iff $S$ admits a more plausible theory than $T$ in describing the concept.

The second desideratum says that theories of sets in the set-theoretical multiverse, ideally, should be ordered with respect to their plausibility, and that two theories are expected to be comparable. One metric, for example, that ensures this desideratum is by looking at the consistency strength of theories to determine their degree of intentionality. If $T$ and $S$ are two `natural' extensions of ZFC, then the degree of intentionality of $T$ and $S$ are comparable. However, Hamkins (2022) claims to have constructed `natural' theories whose consistency strengths are incomparable. He uses a non-standard way of enumerating theories, called {\em cautious enumeration}, and compares diverse alternative cautious enumerations of the ZFC theory, and he argues they are of incomparable consistency strengths to one another. Even if we assume Hamkins is correct, the determination of the degree of intentionality of theories in hierarchical pluralism is not limited to looking at their relative consistency strength. Consistency strength is rather just one metric I gave as an example. Any metric that complies with the listed four desiderata would suffice to have a desirable multiverse structure of set theories for ensuring an equilibrium in pluralist's dilemma. Linearity property is only desired in so far as we want to preserve a robust comparison between, in particular, two competing theories with respect to our best possible theory. Having a partially ordered class of theories does not limit our ability to form a hierarchy of multiverses as we may, for instance, use new additional ordering metrics for further comparison between theories that were originally incomparable with respect to the other metric.
\vspace{0.5cm}

\noindent{\em 3. Coherency.}


(i) (Anti-symmetry) If $T$ and $S$ are two theories with degrees of intentionality, respectively, ${\bf a}$ and ${\bf b}$ such that ${\bf a\leq b}$ and ${\bf b\leq a}$, then ${\bf a=b}$.

(ii) (Transitivity) If $T_1, T_2$, and $T_3$ are theories with degrees of intentionality, respectively, ${\bf a}$, ${\bf b}$, and ${\bf c}$, then whenever ${\bf a}\leq{\bf b}$ and ${\bf b}\leq{\bf c}$, we have that ${\bf a}\leq{\bf c}$.


\vspace{0.5cm}

If we are to take consistency strength for our metric to determine plausibility, there are two special cases which may be, in some sense, thought of as defining the greatest and the least elements in the ordering of the degrees. Any subset of the set of all logical tautologies can be thought to have the {\em greatest} degree of intentionality, whereas any set from which we can prove a contradiction can be thought to have the {\em least} degree of intentionality.\footnote{In fact, these are also the least and the greatest elements of the consistency strength hierarchy.} This is simply due to that any set $T$ of tautologies is already included in our theory and that any model actually satisfies $T$. But we do not necessarily have to include these two degrees as a part of the ordering of multiverse for reasons that will become clear in the next paragraph.

Next is to define truth from the viewpoint of hierarchical pluralism. Suppose we have a conception expressed by a theory. How do we define truth on hierarchical pluralism then? Truth is determined inside the theory with the `maximal' degree of intentionality. For this I shall explain what it means for a theory to have that property. Let $T$ be a theory with degree of intentionality ${\bf a}$. If there exists no known theory $S$ for the same concept with degree of intentionality ${\bf b}$ such that ${\bf a<b}$, then $T$ is said to have a {\em maximal} degree of intentionality. In this case we call $T$ the {\em maximally intended theory} for the conception that is being considered. It should be understood from the first three desiderata that maximally intended theory exists and is unique in the sense if $T$ and $S$ are two maximally intended theories, then the deductive closure of $T$ and $S$ are the same. The maximally intended theory could be the most fruitful theory or the most intrinsically plausible theory that is known to us, depending on the choice of metric and whether we endorse an extrinsic-based naturalistic philosophy or an intrinsic-based position in determining the plausibility of a theory. Essentially it is the best possible theory in hand with respect to the chosen metric. Truth on hierarchical pluralism can be then defined as follows:
\vspace{0.5cm}

\noindent{\em 4. Truth in a hierarchy of theories.}
A statement $\varphi$ is {\em true} iff $\varphi$ is true in the maximally intended theory $T$. If $\psi$ is a statement such that $\psi$ is independent from the maximally intended theory, then there exists a theory $T'$ in which either $\psi$ or $\neg\psi$ is true such that $T'$ has a higher degree of intentionality than that of $T$ and that all statements that are true in $T$ are also true in $T'$.
\vspace{0.5cm}

Due to  G\"{o}del's incompleteness theorem, any recursively axiomatisable theory of sets will be incomplete.\footnote{One might ask what to do with non-recursively enumerable set theories. The axioms of such theories are not even uniformly determined. It is beyond our understanding whether non-recursively enumarable theories are conceivable in the usual sense since we need to go beyond effective methods. Same holds for non-standard theories of sets. For this reason, I consider non-recursively enumerable theories as a part of non-standard theories.} But since there is a fact of the matter whether any given statement is true or false, even so that it might be undecidable in the maximally intended model, it does not necessarily mean that it is absolutely undecidable. So hierarchical pluralism rejects {\em absolute} undecidability, for it is assumed there is an extension of the maximally intended theory in which the undecidable statement is settled. This is another reason that hierarchical pluralism is a realist based multiverse theory in the truth value sense. But it is still congruent with the notion of independence in the temporal sense. This makes hierarchical pluralism compatible with the naturalistic approach to set-theoretical practice.

These four desiderata are not imposed on a specific metric such as, for instance, relative consistency strength. These are rather a set of requirements that one would expect from a hierarchy of different set theoretical universes while preserving epistemological realism, i.e., realism in truth value, as much as possible. Due to pluralist's dilemma, however, it does not seem possible to fully satisfy either of the extremes independent of the metric being considered. Hierarchical pluralism offers an equilibrium between the two extremes---objectivity of truth vs. richness of the ontology---by allowing a multiverse view to ensure abundance, while accepting a degree of intentionality for every theory in the multiverse for preserving epistemological realism with the aid of defining truth relative to the maximally intended theory.

Hierarchical pluralism for sets could either be conceptualized based on `intrinsic plausibility' in the sense of intrinsic features of models or theories and how it conforms to the full conception of sets, or it could be conceptualized based on the `fruitfulness' of theories. In any case, hierarchical pluralism is a view where we have a series of successive grasps towards the true conception regardless of whether this approximation relies more on extrinsic or intrinsic properties of theories. In Friend's `growing foundation' sense, maximally intended theory can be outgrown or replaced as we discover more mathematical results, particularly when we select theories by their extrinsic properties. Imagine the following example. Suppose that ZFC were our current maximally intended theory about the concept of `set'. If we ever select, by practice or as a result of an intuitive insight, $\textrm{ZF}+\textrm{AD}$ over ZFC, then sentences that are entailed from $\textrm{ZF}+\textrm{AD}$ will be {\em preferred} over the sentences entailed from ZFC.\footnote{AD stands for the {\em Axiom of Determinacy}, which is roughly the statement that for any given two-player infinite game, there exists a winning strategy for one of the players. There is a conflict between AD and the Axiom of Choice. Now both have defects and merits. One of the most famous paradoxical results of the Axiom of Choice is the {\em Banach-Tarski Theorem}: For any given ball in the 3-dimensional space, there exists a decomposition of the ball into a finite number of disjoint subsets such that when put back together in a particular way, yields two identical copies of the original ball. However, AD is not entirely paradox-free either. It follows from AD that a set $A$ with cardinality $2^{\aleph_0}$ can be partitioned into its disjoint subsets such that the number of elements in the partition is greater than the number of  elements of $A$. See Halbeisen and Shelah \cite{HalbeisenShelah2001} (2001), p. 258, for more details.} In this case, we say that the previous maximally intended theory, i.e. ZFC, is {\em surpassed} by $\textrm{ZF}+\textrm{AD}$. So ZFC is replaced now by $\textrm{ZF}+\textrm{AD}$ which becomes our new maximally intended theory. It may be the case for later that we decide to accept a higher axiom of infinity $\Phi$ by, say, naturalistic reasons such that $\Phi\rightarrow\neg\textrm{AD}$. Then, $\textrm{ZF}+\Phi$ will have a higher degree of intentionality than that of $\textrm{ZF}+\textrm{AD}$, and our previous insights about AD turns out implausible at this stage. In general, suppose that $\textrm{ZFC}+\varphi$ is preferred over $\textrm{ZFC}+\neg\varphi$ for pragmatic reasons that the degree of intentionality of $\textrm{ZFC}+\varphi$ is higher than that of $\textrm{ZFC}+\neg\varphi$. If we find some $\psi$ such that $\psi\rightarrow\neg\varphi$ and that the degree of intentionality of $\textrm{ZFC}+\psi$ is greater than that of $\textrm{ZFC}+\varphi$, then $\textrm{ZFC}+\neg\varphi$ is now preferred over $\textrm{ZFC}+\varphi$. Adding new axioms merely based on extrinsic features does not necessarily mean that the new maximally intended theory {\em really} has a higher degree of intentionality. When we rely on extrinsic justifications more, our selections may later turn out to be made poorly. Suppose we decide to add a large cardinal axiom $\Psi$ in ZFC and we later find another axiom candidate which falsifies $\Psi$. This means that adding $\Psi$ in the first place was a poor selection. So taking into account, majority of the time, extrinsic features of theories for adding new axioms leads to a more dynamic picture in the multiverse hierarchy, and hence our theories may get replaced by new selections. Such replacements can occur in both intrinsic-based and extrinsic-based hierarchical pluralism. The difference is that when we take the extrinsic route, the plausibility ordering of theories becomes more dependent on our pragmatic intentions and the mathematical practice. Whereas on the intrinsic route, whenever a new theory of sets is proposed, say due to some discovery of an obviously true statement, the place of the new theory in the ordering is fixed as if it is already {\em pre-determined} in advance independent of our pragmatic preferences. This is due to the fact that since the full conception exists on hierarchical pluralism, there really {\em exists} an ordering of universes with respect to their intrinsic plausibility.

Let me now go back to Balaguer's treatment of pluralism and explain how hierarchical pluralism extends it and offers an alternative solution to find an equilibrium in pluralist's dilemma. The kind of Platonism that Balaguer suggests for Platonists to endorse is what he calls {\em intention-based Platonism} (IBP) for which he defines truth in the following manner:

\begin{changemargin}{1cm}{1cm}
{\em IBP $($short for ``{\em intention-based platonism}"$)$: A pure mathematical sentence $S$ is {\em true} iff it is true in {\em all} the parts of the mathematical realm that count as intended in the given branch of mathematics $($and there is at least one such part of the mathematical realm$)$; and $S$ is {\em false} iff it is false in all such parts of the mathematical realm $($or there is no such part of the mathematical realm$)$; and if $S$ is true in some intended parts of the mathematical realm and false in others, then there is no fact of the matter whether it is true or false.\footnote{Balaguer, pp. 145-146, 2009.}}
\end{changemargin}

It is in fact this type of Platonism which Balaguer claims is consistent with mathematical pluralism. Regarding IBP, he concludes the following: 

\begin{changemargin}{1cm}{1cm}
{\em Given that IBP is a {\em plenitudinous} version of platonism, we get the following result: which $($consistent$)$ mathematical sentences are true is wholly determined by our intentions. This is because $($a$)$ mathematical truth comes down to truth in intended structures, and $($b$)$ which mathematical structures count as intended is wholly determined by our intentions.\footnote{Balaguer, p. 393, 2017.}}
\end{changemargin}
 
The justification used in part (b) seems to aim for a naturalistic philosophy if what is meant by `intentions' is `fruitfulness'. However, if it the word `intention' refers to `mathematical intuition', then we need to presuppose the existence of what is referred to as the full conception of the theory.

Hierarchical pluralism settles independent statements that are claimed to have `no fact of the matter whether they are true or false' in Balaguer's IBP. This was provided in the fourth desideratum. So this kind of realist-based pluralism extends Balaguer's IBP in the sense that it is granted that there exists {\em some} theory in which the truth value of the given independent statement is settled.\footnote{Woodin's axioms for the $V=\textrm{Ultimate-}L$ is a great example to this, where he claims that all {\em known} large cardinal axioms that were known to be independent are settled in the ultimate extension of $L$. For details, see Woodin (2017).}


Balaguer's claim is that a certain version of Platonism and anti-Platonism both entail pluralism, and so he concludes that pluralism must be true.\footnote{See Balaguer (2017), pp. 392-393, for more details.} For Balaguer, if $M_1$ and $M_2$ are two models of set theory---possibly `equally' consistent with our full conception of sets---and if $M_1\models \varphi$ and $M_2\models\neg\varphi$, then there is no fact of the matter as to whether $\varphi$ is true or false. Again, this relies on the equal-status interpretation of the multiverse theory. So on the IBP interpretation, the status of CH is somewhat left open since set theorists have experiences in CH and $\neg\text{CH}$ universes, both of which can be claimed as intended models of set theory. In fact on IBP, in some cases, there {\em really} is no fact of the matter whether CH is true or false. But I think this does not convey the most basic spirit of {\em epistemological realism}, the idea that every mathematical statement really {\em has} a truth value independent of our sensory experience, mind, language, etc. One may worry about that IBP does not give us a robust kind of objectivity. Balaguer counters this worry, as we quoted above, by claiming that consistent mathematical structures are determined entirely by our intentions. But then on IBP, mathematical reality is bound by our intentions (hence, our mind and language) about the mathematical concepts. It seems to me that IBP offers a type of {\em language-dependent} realism. On hierarchical pluralism, ideally, either CH or $\neg\text{CH}$ universe has a higher degree of intentionality. Thus, one of them is the preferred universe over the other. We may not know in reality which one is it, but it is granted in a non-uniform way that there at least exists a structure in which CH is a settled proposition. So based on this schema, every statement {\em is} settled in one way or another. It appears to me that this hierarchical schema is a reasonable way of preserving epistemological realism (and the mathematical objectivity) within the multiverse theory of sets. Hierarchical pluralism serves as a reconcilation of epistemological realism, i.e. realism in truth value, and the multiverse theory. In fact, Balaguer's IBP accomplishes this to a certain extent. However, he also claims that the truth value of CH, or any other independent statement of similar type, `has no fact of the matter as to whether it is true or not'. Hierarchical pluralism fills this gap by looking at the truth value in the best possible theory, or in other words, in the maximally intended theory.

On hierarchical pluralism, then, we order the theories with respect to their intendedness and define truth instead in the most preferred theory. This amounts to some objectivity, not with respect to an arbitrary universe, but with respect to the most plausible universe in this sense. Two competing theories may be both plausible, but one of them, ideally, is more preferred than the other. Statements that are satisfied in two competing models, likewise, may be both true respectively in those models, but the statements that are true in the preferred model have a higher intendedness than the other. In other words, the credibility of truth should vary in degrees in proportional to the intentionality of theories. This is the basic idea of hierarchical pluralism. The set theorist, in this way, settles the pluralist's dilemma by balancing out realism in truth value with the aid of defining truth in maximally intended theory and abundancy of mathematical ontology by means of allowing the existence of other theories with lower degrees of intentionality. 

One might ask whether multiverse related concerns affect geometry the way it does set theory. Perhaps we can ask the following question: Is there a pluralist's dilemma for geometry? The debate about geometry vs. set theory, as to why a pluralistic view on geometry is more different than in set theory, is an interesting problem that, I believe, deserves an independent study in its own right. As I have mentioned earlier, I think the conceptual line between Euclidean and non-Euclidean geometries is more firm than that of between different set theories. The reason for this is the following. We tend to assume that independent statements of geometry, particularly the parallel postulate for instance, are actually intended to be posited in the context of Euclidean conception of space. We cannot be sure in what set context CH was formulated. Two space conceptions, i.e., Euclidean and non-Euclidean, seem entirely different to us through satisfied axioms, but most importantly it is more straightforward to distinguish the `space' conception of one from the other. This is not the case for sets due to that we don't have spatial or visual tools for differentiating set concepts of the CH and non-CH universes. The reason is that I think set theory is too complicated to have a `natural' structure and hence it gets harder to compare with other competing theories. We automatically assume the Euclidean space conception when we talk about the parallel postulate, but we don't assume any natural set conception---at least not as much as we do so in geometry---when we talk about CH.


{\small 
\begin{thebibliography}{99}

\bibitem{AntosFriedman}
Carolin Antos, Sy David Friedman, Radek Honzik, Claudio Ternullo, {\em Multiverse Conceptions in Set Theory}, Synthese, {\bf 192}(8), pp. 2463-2488 (2015).

\bibitem{ArrigoniFriedman}
Tatiana Arrigoni, Sy David Friedman, {\em The Hyperuniverse Programme}, Bull. Symbolic Logic, {\bf 19}(1), pp. 77-96 (2013).
 
\bibitem{Balaguer1998}
Mark Balaguer, {\bf Platonism and anti-platonism in mathematics}, Oxford University Press, New York (1998). 
 
\bibitem{Balaguer2009}
Mark Balaguer, {\em Fictionalism, Theft, and the Story of Mathematics}, Philosophia Mathematica, {\bf 17}(2), pp. 131-162 (2009). 
 
\bibitem{Balaguer2016}
Mark Balaguer, {\em Full-blooded platonism}; in {\bf An historical introduction to the philosophy of mathematics}, R. Marcus and M. McEvoy (eds.),  pp. 719-732, Bloomsbury Publishing, New York (2016). 
 
\bibitem{BalaguerPluralism}
Mark Balaguer, {\em Mathematical Pluralism and Platonism}, J. Indian Counc. Philos. Res., {\bf 34}(2), pp. 379-198, Springer (2017).

\bibitem{Barton}
Neil Barton, {\em Multiversism and Concepts of Set: How much Relativism is acceptable?}; in {\bf Objectivity, Realism, and Proof}, Francesca Boccuni \& Andrea Sereni (eds.), pp. 189-209, Springer (2016).

\bibitem{BenacerrafPutnam}
Paul Benacerraf, Hilary Putnam (ed.), {\bf Philosophy of Mathematics: Selected readings}, Cambridge University Press, New York (1983).

\bibitem{Berry2022}
Sharon Berry, {\bf A Logical Foundation for Potentialist Set Theory}, Cambridge University Press, New York (2022).

\bibitem{Bishop}
Eric Bishop, {\bf Foundations of Constructive Analysis}, McGraw-Hill, New York (1967).

\bibitem{Button}
Tim Button, Sean Walsh, {\bf Philosophy and Model Theory}, Oxford University Press, Oxford (2018).

\bibitem{Caicedo2011}
Andr\'{e}s E. Caicedo. {\em (Non?)-linearity of the consistency strength ordering in ZF}, MathOverflow answer, \url{https://mathoverflow.net/q/59800}, (2011).

\bibitem{Cohen}
Paul Cohen, {\em The Independence of the Continuum Hypothesis, I, II}, Proc. Nat. Acad. Sci. U.S.A., {\bf 50}, pp. 1143-1148 (1963). 

\bibitem{DaviesPluralism}
Brian Davies, {\em Pluralism in Mathematics}, Philos. Trans. A. Math. Phys. Eng. Sci. {\bf 363}(1835), pp. 2449-60 (2005).

\bibitem{FriendPluralism}
Mich\`{e}le Friend, {\bf Pluralism in Mathematics: A New Position in Philosophy of Mathematics}, Springer (2014).

\bibitem{GodelCH}
Kurt G\"{o}del, {\em Consistency-proof for the Generalized Continuum Hypothesis}, Proc. Nat. Acad. Sci. U.S.A., {\bf 25}, pp. 220-224 (1939).

\bibitem{Godel1964}
Kurt G\"{o}del, {\em What is Cantor's Continuum Problem?}, 1964; reprinted in {\bf Philosophy of Mathematics Selected Readings}, Paul Benacerraf and Hilary Putnam (eds.), pp. 470-485, Cambridge University Press (1983).

\bibitem{HalbeisenShelah2001}
Lorenz Halbeisen, Saharon Shelah, {\em Relations Between Some Cardinals in the Absence of the Axiom of Choice}, Bulletin of Symbolic Logic, {\bf 7}(2), pp. 237-261 (2001).

\bibitem{HamkinsMultiverse}
Joel David Hamkins, {\em The Set-theoretic Multiverse}, Review of Symbolic Logic, {\bf 5}, pp. 416-449 (2012).

\bibitem{Hamkins2022}
Joel David Hamkins, {\em Nonlinearity and illfoundedness in the hierarchy of large cardinal consistency strength}, available at \url{https://arxiv.org/abs/2208.07445} (2022).

\bibitem{HamkinsLinnebo2019}
Joel David Hamkins, {\O}ystein Linnebo, {\em The modal logic of set-theoretic potentialism and the potentialist maximality principles}, Review of Symbolic Logic, {\bf 15}(1), pp. 1-35 (2022).

\bibitem{HellmanPluralism}
Geoffrey Hellman, John Lane Bell, {\em Pluralism and the Foundations of Mathematics}; in {\bf Minnesota Studies in the Philosophy of Science Series}, John Lane Bell (ed.), Vol. 19, pp. 64-79, University of Minnesota Press (2006).

\bibitem{Jensen}
Ronald Jensen, {\em Inner Models and Large Cardinals}, The Bulletin of Symbolic Logic, {\bf 1}(4), pp. 393-407 (1995).


\bibitem{Koellner2011}
Peter Koellner, {\em Independence and Large Cardinals}, The Stanford
Encyclopedia of Philosophy, Ed. by Edward N. Zalta (2011).

\bibitem{KoellnerHamkinsMultiverse}
Peter Koellner, {\em Hamkins on the Multiverse}, unpublished notes, available at \url{http://logic.harvard.edu/EFI_Hamkins_Comments.pdf} (2013).

\bibitem{Linnebo2013}
{\O}ystein Linnebo, {\em The potential hierarchy of sets}, Review of Symbolic Logic, {\bf 6}(2), pp. 205-228 (2013).

\bibitem{Maddy1997}
Penelope Maddy, {\bf Naturalism in Mathematics}, Oxford University Press (1997).


\bibitem{Maddy2011DefendingAxioms}
Penelope Maddy, {\bf Defending the Axioms: On the Philosophical Foundations of Set Theory}, Oxford University Press, New York (2011).

\bibitem{Maddy2017}
Penelope Maddy, {\em Set-theoretic foundations}, Contemporary Mathematics, {\bf 690}, pp. 289-322 (2017).

\bibitem{Magidor2012}
Menachem Magidor, {\em Some set theories are more equal}, unpublished notes, available at \url{http://logic.harvard.edu/EFI_Magidor.pdf} (2012).

\bibitem{Mirimanoff1917}
Dmitry Mirimanoff, {\em Les antinomies de Russell et de Burali-Forti et le probleme fondamental de la theorie des ensembles}, L'Enseignement Math\'{e}matique, {\bf 19}, pp. 37-52 (1917).


\bibitem{Poincare1902}
Henri Poincar\'{e}, {\bf La Science et l'Hypoth\`{e}se}, 1902; trans. by William John Greenstreet, The Walter Scott Publishing, Co., Ltd, New York (1905).

\bibitem{PriestMathPluralism}
Graham Priest, {\em Mathematical Pluralism}, Logic Journal of the IGPL, {\bf 21}(1), pp. 4-13 (2013).


\bibitem{Shapiro2000}
Stewart Shapiro, {\bf Thinking about mathematics: Philosophy of mathematics}, Oxford University Press, New York (2000).

\bibitem{Shelah}
Saharon Shelah, {\em Logical Dreams}, Bullettin of the American Mathematical Society, {\bf 40}(2), pp. 203-228 (2003). 

\bibitem{Simpson2009}
Stephen G. Simpson, {\em The G\"{o}del Hierarchy and Reverse Mathematics}, available at \url{http://www.personal.psu.edu/t20/papers/gh/} (2009).

\bibitem{SteelInnerModel}
John Steel, {\em An outline of inner model theory}; in {\bf Handbook of Set Theory}, Akihiro Kanamori, Matthew Foreman (eds.), vol. 3, pp. 1595-1684, Springer, Berlin (2010).

\bibitem{Steel2014}
John Steel, {\em G\"{o}del's program}, in {\bf Interpreting G\"{o}del}, Juliette Kennedy (ed.), Cambridge University Press, pp. 153-179 (2014).

\bibitem{Ternullo2019}
Claudio Ternullo, {\em Maddy On The Multiverse}, in {\bf Reflections on the Foundations of Mathematics}, Deniz Sarikaya, Deborah Kant \& Stefania Centrone (eds.), pp. 43-78, Springer Verlag (2019).

\bibitem{WoodinUltimateL}
William Hugh Woodin, {\em In Search of Ultimate-$L$: The 19th Midrasha Mathematicae Lectures}, The Bulletin of Symbolic Logic, {\bf 23}(1), pp. 1-109 (2017).

\bibitem{Wright1986} 
Crispin Wright, {\bf Realism, meaning and truth}, Blackwell, Oxford (1986).

\bibitem{Zermelo}
Ernst Zermelo, {\em \"{U}ber Grenzzahlen und Mengenbereiche: neue Untersuchungen \"{u}ber die Grundlagen der Mengenlehre}, Fundamenta Mathematicae, {\bf 16}, pp. 29-47 (1930). 


%
%
%
%
%
%
%
%
%
%
%
%
%
%
%
%
%
%
%
%
%
%
%
%
%
%
%
%
%
%
%
%
%
%
%
%

%
%
%
%




\end{thebibliography}
}
\end{document}